\documentclass[12pt]{amsart}
\usepackage{amscd,amssymb}
\usepackage{upref}
\usepackage{pstcol,pst-3d}

\psset{unit=0.7cm,linewidth=0.4pt,arrowsize=1.5pt 4}
\def\vertex{\pscircle[fillstyle=solid,fillcolor=black]{0.0714}}
\newpsstyle{fatline}{linewidth=1pt}
\newpsstyle{fyp}{fillstyle=solid,fillcolor=verylight}
\definecolor{verylight}{gray}{0.95}
\definecolor{light}{gray}{0.9}
\definecolor{medium}{gray}{0.85}

\let\oldlabel=\label
\def\prellabel{\marginparsep=1em
    \def\label##1{\oldlabel{##1}\ifmmode\else\ifinner\else
         \marginpar{{\footnotesize\ \\ \tt
                    ##1}}\fi\fi}}

\def\bF{{\textbf{\textup{F}}}}
\def\bG{{\textbf{\textup{G}}}}

\let\frak=\mathfrak
\let\cal\mathcal
\let\Bbb\mathbb

\let\la=\langle
\let\ra=\rangle

\def\cd#1{$\textup{CD}_{#1}$}

\def\ggr.aut{{\operatorname{gr.aut}}}

\def\L{\operatorname{L}}
\def\SL{\operatorname{SL}}
\def\int{{\operatorname{int}}}
\def\inv{{\operatorname{inv}}}
\def\vert{{\operatorname{vert}}}
\def\conv{\operatorname{conv}}
\def\het{\operatorname{ht}}
\def\E{\operatorname{E}}
\def\V{\operatorname{V}}

\def\St{\operatorname{{\Bbb S}t}}

\def\path{\operatorname{path}}
\def\B{\operatorname{B}}

\def\T{\operatorname{\bf T}}

\def\Aff{\operatorname{aff}}
\def\End{\operatorname{End}}
\def\gp{\operatorname{gp}}

\def\Ker{{\operatorname{Ker}}}
\def\Coker{{\operatorname{Coker}}}

\def\Hom{\operatorname{Hom}}
\def\Pol{\operatorname{Pol}}
\def\GL{\operatorname{GL}}
\def\Qu{\operatorname{Q}}
\def\Col{\operatorname{Col}}

\def\Q{\qedsymbol\kern1pt}
\def\sqq#1#2{{\hbox{\kern 0.5pt\vbox{\hbox{\kern#1ex\vrule width #2pt height #1ex}%
     \hrule height #2pt}\kern 1.5pt}}}
\def\sq{{\mathchoice{\sqq1{0.5}}{\sqq1{0.5}}{\sqq{0.9}{0.4}}
                    {\sqq{0.7}{0.3}}}}

\def\vt{{\kern1pt|}}
\def\RR{{\mathbb R}}

\def\ZZ{{\mathbb Z}}
\def\NN{{\mathbb N}}
\def\TT{{\mathbb T}}
\def\AA{{\mathbb A}}
\def\EE{{\mathbb E}}

\def\vv{{\mathbb V}}
\def\fP{{\mathfrak P}}
\def\cG{{\mathcal G}}

\let\epsilon=\varepsilon
\let\phi=\varphi
\let\theta=\vartheta

\advance\headheight1.15pt

\newtheorem{lemma}{Lemma}[section]

\newtheorem{theorem}[lemma]{Theorem}
\newtheorem{conjecture}[lemma]{Conjecture}
\newtheorem{proposition}[lemma]{Proposition}

\theoremstyle{definition}
\newtheorem{definition}[lemma]{Definition}
\newtheorem{remark}[lemma]{Remark}

\begin{document}

\title[Polytopes \and $K$-theory]
{Polytopes and $K$-theory}

\author{Winfried Bruns \and Joseph Gubeladze}

\address{Universit\"at Osnabr\"uck,
FB Mathematik/Informatik, 49069 Osnabr\"uck, Germany}
\email{Winfried.Bruns@mathematik.uni-osnabrueck.de}

\address{A. Razmadze Mathematical Institute, Georgian Academy of Sciences,
1 Alexidze str., Tbilisi 0193, Republic of Georgia\ \&\newline
Department of Mathematics, San Francisco State University, San
Francisco, CA 94132, USA} \email{soso@math.sfsu.edu}

\thanks{Both authors were supported by Mathematical
Sciences Research Institute}


\maketitle

\hfill\underline{\emph{\small To the memory of Professor George
Chogoshvili (1914-1998)}}

\

\

\begin{flushright}
\emph{\tiny Every so often you should try a damn-fool experiment\ \ ---  }\\
{\tiny from J. Littlewood's} \textsc{\tiny A Mathematician's
Miscellany}
\end{flushright}

\section{Introduction}\label{INTR}

We overview results from our experiment
of merging two seemingly unrelated disciplines --
higher algebraic $K$-theory of rings and the theory of lattice
polytopes. The usual $K$-theory is the ``theory of a unit
simplex''.

The text is based on the works \cite{BrG1,BrG5,BrG6}.

At the end of the paper we propose a general conjecture on the
structure of higher polyhedral $K$-groups for certain class of
polytopes for which the coincidence of Quillen's and Volodin's
theories is known.

\emph{All rings, considered below, are commutative and for a ring
$R$ its multiplicative group of units is denoted by $R^*$.}

\section{Motivation and applications}\label{motivations}

To defuse the impression on the experiment to be too damn-fool,
here we describe the motivation behind our polyhedral $K$-theory.

Demazure's paper \cite{D} that initiated the theory of toric
varieties in the early 1970s gave an exhaustive description of the
automorphism group of a complete smooth toric variety. (Much later
this was extended to arbitrary complete toric varieties by Cox
\cite{C} and Buehler \cite{Bu}.) Theorem \ref{PLg} below gives an
analogous result for the graded automorphism group of the affine
cone over a projective toric variety, not necessarily smooth. As
explained in Section \ref{PLG}, this approach leads to
\emph{polytopal generalizations} of the groups $\GL_n(k)$, $k$ a
field, and the standard fact that $\SL_n(k)=\E_n(k)$. Our
motivating question is: to what extent the \emph{polytopal linear
groups} and the associated higher $K$-groups resemble the ordinary
$K$-groups? We work with the techniques of Quillen's +
construction and Volodin's definition of higher $K$-groups. This
seems the only possible framework in our essentially non-additive
situation.

On the level of $K_2$, polyhedral $K$-theory can be thought of as
complementary to the theory of \emph{universal Chevalley groups}
\cite{KStn,Stb,Stn}. This is so because polytopal linear groups
are semidirect products of unipotent groups and reductive groups
of type $A_l$, see \cite[Section 1]{BrG5}.

For higher groups one is naturally led to the study of the
integral homology of interesting examples of linear groups, see
Section \ref{POLYG}.

As an application to toric geometry, we have obtained results on
\emph{retractions} of toric varieties \cite{BrG2}, automorphisms
of \emph{arrangements} of toric varieties \cite{BrG3}, and
\emph{autoequivalences} of the category of toric varieties
\cite{BrG4}.

\section{Polytopes, their algebras, and their linear groups}\label{POLYTOPES}

\subsection{General polytopes}\label{GENPOL}

By a \emph{polytope} $P\subset\RR^n$, $n\in\NN$, we always mean a
\emph{finite convex} polytope, i.~e.\ $P$ is the convex hull of a
finite subset $\{x_1,\ldots,x_k\}\subset\RR^n$:
$$
P=\conv(x_1,\ldots,x_k):=\{a_1x_1+\dots+a_kx_k\ :\ 0\le
a_1,\dots,a_k\le1,\\ a_1+\cdots+a_k=1\}.
$$

Polytopes of dimension 1 are called \emph{segments} and those of
dimension 2 are called \emph{polygons}.

The \emph{affine hull} $\Aff(X)$ of a subset $X\subset\RR^n$ is
the smallest affine subspace of $\RR^n$ containing $X$. If $\dim
\Aff(X)=k-1$ for a subset $X=\{x_1,\dots,x_k\}$ of cardinality
$k$, then $x_1,\dots,x_k$ are \emph{affinely independent} and the
polytope $P=\conv(x_1,\dots,x_k)$ is called a \emph{simplex}.

For a halfspace $\mathcal H\subset\RR^n$ containing $P$, the
intersection $P\cap\partial\mathcal H$ of $P$ with the affine
hyperplane $\partial\mathcal H$ bounding $\mathcal H$ is called a
\emph{face} of $P$. The polytope itself is also considered as a
face.

The faces of $P$ are themselves polytopes. Faces of dimension $0$
are \emph{vertices} and those of codimension $1$ (i.~e.\ of
dimension $\dim P-1$) are called \emph{facets}. A polytope is the
convex hull of the set $\vert(P)$ of its vertices. If $\dim
P\subset\RR^n$ has dimension $n$, then there is a unique halfspace
$\mathcal H$ for each facet $F\subset P$ such that
$P\subset\mathcal H$ and $\partial\mathcal H\cap P=F$.

\subsection{Lattice polytopes}\label{LATPOL} A polytope $P\subset\RR^n$ is
called a \emph{lattice polytope} if the vertices of $P$ belong to
the integral lattice $\ZZ^n$. More generally, a lattice in $\RR^n$
is a subset $\cG=x_0+\cG_0$ with $x_0\in\RR^n$ and an  additive
subgroup $\cG_0$ generated by $n$ linearly independent vectors. A
polytope $P$ with $\vert(P)\subset \cG$ is called a $\cG$-polytope
if the vertices of $P$ belong to $\cG$. However, since all the
properties of $\cG$-polytopes we are interested in remain
invariant under an affine automorphism of $\RR^n$ mapping $\cG$ to
$\ZZ^n$, we can always assume that our polytopes have vertices in
$\ZZ^n$. More generally, lattice polytopes $P$ and $Q$ that are
isomorphic under an integral-affine equivalence of $\Aff(P)$ and
$\Aff(Q)$ are equivalent objects in our theory. We simply speak of
\emph{integral-affinely equivalent polytopes}.

Faces of a lattice polytope are again lattice polytopes.

For a lattice polytope $P\subset\RR^n$ we put $\L_P=P\cap\ZZ^n$. A
simplex $\Delta$  is called \emph{unimodular} if
$\sum_{z\in\vert(\Delta)}\ZZ(z-z_0)$ is a direct summand of
$\ZZ^n$ for some (equivalently, every) vertex $z_0$ of $\Delta$.
All unimodular simplices of dimension $n$ are integral-affinely
equivalent. Such a simplex is denoted by $\Delta_n$ and called a
\emph{unit $n$-simplex}. Standard realizations of $\Delta_n$ are
$\conv(O,e_1,\dots,e_n)\subset\RR^n$ or
$\conv(e_1,\dots,e_{n+1})\subset\RR^{n+1}$. ($e_i$ is the $i$th
unit vector.)

There is no loss of generality in assuming that a given lattice
polytope $P$ is full dimensional (i.~e.\ $\dim P=n$) and that
$\ZZ^n$ is the smallest affine lattice containing $\L_P$.  In
fact, we choose $\Aff(P)$ as the space in which $P$ is embedded
and fix a point $x_0\in\L_P$ as the origin. Then the lattice
$x_0+\sum_{x\in\L_P}\ZZ(x-x_0)$ can be identified with $\ZZ^r$,
$r=\dim P$.

Under this assumption let $F$ be a facet of $P$ and choose a point
$z_0\in F$. Then the subgroup
$$
F_\ZZ:=(-z_0+\Aff(F))\cap\ZZ^n\subset\ZZ^n
$$ is isomorphic to $\ZZ^{n-1}$. Moreover, there is a unique group
homomorphism $\la F,-\ra:\ZZ^n\to\ZZ$, written as $x\mapsto\la
F,x\ra$, such that $\Ker(\la F,-\ra)=F_\ZZ$, $\Coker(\la
F,-\ra)=0$, and on the set $\L_P$, $\la F,-\ra$ attains its
minimum $b_F$  at the lattice points of $F$.

The $\ZZ$-linear form $\la F,-\ra$ can be extended in a unique way
to a linear function on $\RR^n$. The description of $P$ as an
intersection of halfspaces yields that $x\in P$ if and only if
$\la F,x\ra\ge b_F$ for all facets $F$ of $P$.

\emph{All polytopes, considered below, are lattice polytopes.}

\subsection{Column structures}\label{COLSTR} 
\begin{figure}[htb]
\begin{center}
\begin{pspicture}(0,0)(5,3)
\pspolygon[style=fyp](0,0)(5,0)(5,2)(4,3)(2,3)(1,2)
\psline[linecolor=gray](2,0)(2,3)
\psline[linecolor=gray](3,0)(3,3)
\psline[linecolor=gray](4,0)(4,3)
\psline[linecolor=gray](1,0)(1,2) \multirput(0,0)(1,0){6}{\vertex}
\multirput(1,1)(1,0){5}{\vertex} \multirput(1,2)(1,0){5}{\vertex}
\multirput(2,3)(1,0){3}{\vertex} \rput(1.3,0.5){$v$}
\psline[style=fatline]{->}(1,1)(1,0.05)
\end{pspicture}
\end{center}
\caption{A column structure}\label{FigCol}
\end{figure}
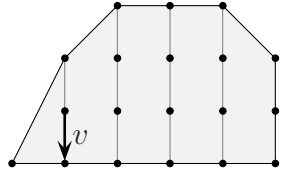

Let $P\subset\RR^n$ be a
polytope. A nonzero element $v\in\ZZ^n$ is called a \emph{a column
vector} for $P$ if there exists a facet $F\subset P$ such that
$x+v\in P$ whenever $x\in\L_P\setminus F$. In this situation $F$
is uniquely determined and called the \emph{base facet} of $v$. We
use the notation $F=P_v$. The set of column vectors of $P$ is
denoted by $\Col(P)$. A \emph{column structure} is a pair of type
$(P,v)$, $v\in\Col(P)$. Figure \ref{FigCol} gives an example of a
column structure.
Familiar examples of column structures are the unit simplices
$\Delta_n$ with their edge vectors.

\subsection{Polytopal semigroups and their rings}\label{POLRINGS}

To a polytope $P\subset\RR^n$ one associates the additive
subsemigroup $S_P\subset\ZZ^{n+1}$, generated by $\{(z,1)\ :\
z\in\L_P\}\subset\ZZ^{n+1}$. Let $C_P\subset\RR^{n+1}$ be the cone
$\{az\ :\ a\in\RR_+,\ z\in P\}$. Then $C_P$ is the convex hull of
$S_P$. It is a \emph{finite rational pointed} cone. In other
words, $C_P$ is the intersection of a finite system of halfspaces
in $\RR^{n+1}$ whose boundaries are rational hyperplanes
containing the origin $O\in\RR^{n+1}$, and there is no affine line
contained in $C_P$.

As in Subsection \ref{LATPOL}, there is no loss of generality in
assuming that $\ZZ^n$ is the lattice spanned affinely by $\L_P$ in
$\RR^n$. This is equivalent to $\gp(S_P)=\ZZ^{n+1}$.

While the points $x\in \L_P$ are identified with
$(x,1)\in\ZZ^{n+1}$, a column vector $v$ is to be identified with
$(v,0)\in\ZZ^{n+1}$.

Let $F$ be a facet of $P$. We use the function $\la F,-\ra$ to
define the \emph{height} of
$x=(x',x'')\in\RR^{n+1}=\RR^n\times\RR$ above the hyperplane
$\mathcal H$ through the facet $C_F$ of $C_P$ by setting
$$
\het_F(x)=\la F,x'\ra-x''b_F.
$$
For lattice points $x$ the function $\het_F$ counts the number of
hyperplanes between $\mathcal H$ and $x$ (in the direction of $P$)
that are parallel to, but different from $\mathcal H$ and pass
through lattice points. If $v$ is a column vector, then $\het_v$
stands for $\het_{P_v}$. Moreover, we are justified in calling
$\het_F(v,0)=\la F,v\ra$ the \emph{height of $v$ with respect to
$F$}, since $v$ is identified with $(v,0)$.

Although the semigroup $S_P$ may miss some integral points in the
cone $C_P$ this cannot happen on the segments parallel to a column
vector $v$. More precisely, the following holds:
\begin{equation}\label{COLUMNS}
z+v\in S_P\ \text{for all}\ z\in S_P\setminus C_{P_v}.
\end{equation}
($C_{P_v}\subset C_P$ is the face subcone, corresponding to
$P_v$.)

Let $R$ be a ring and $P\subset\RR^n$ a lattice polytope. The
semigroup ring $R[P]:=R[S_P]$ -- the \emph{polytopal $R$-algebra
of $P$} -- carries a graded structure $R[P]=R\oplus
R_1\oplus\cdots$ in which $\deg(x)=1$ for all $x\in\L_P$. By
definition of $S_P$ it follows that $R_1$ generates $R[P]$ over
$R$.

We are interested in the group $\ggr.aut_R(P)$ of graded
$R$-algebra automorphisms of $R[P]$. For a field $R=k$ the group
$\ggr.aut_k(P)$ is naturally a $k$-linear group. In fact, it is a
closed subgroup of $\GL_m(k)$, $m=\#\L_P$. We call $\ggr.aut_k(P)$
the \emph{polytopal $k$-linear group} of $P$. Its structure will
be given in Theorem \ref{PLg}.

In the special case when $P$ is a unimodular simplex, the ring
$R[P]$ is isomorphic to a polynomial algebra $R[X_1,\ldots,X_m]$,
$m=\#\L_P$. Therefore, the category $\Pol(R)$ of polytopal
$R$-algebras and graded homomorphisms between them contains a full
subcategory that is equivalent to the category of free
$R$-modules.

\subsection{Polytopal linear groups}\label{PLG}

Assume $R$ is a ring and $P$ a polytope. Let $(P,v)$ be a column
structure and $\lambda\in R$. As pointed out above, we identify
the vector $v$ with the degree $0$ element $(v,0)\in\ZZ^{n+1}$,
and further with the corresponding monomial in $R[\ZZ^{n+1}]$.
Then we define a mapping from $S_P$ to $R[\ZZ^{n+1}]$ by the
assignment
$$
x\mapsto (1+\lambda v)^{\het_v x}x.
$$
Since $\het_v$ is a group homomorphism $\ZZ^{n+1}\to\ZZ$, our
mapping is a homomorphism from $S_P$ to the multiplicative monoid
of $R[\ZZ^{n+1}]$. Now it is immediate from (\ref{COLUMNS}) in
Subsection \ref{POLRINGS} that the (isomorphic) image of $S_P$
lies actually in $R[P]$. Hence this mapping gives rise to a graded
$R$-algebra endomorphism $e_v^\lambda$ of $R[P]$ preserving the
degree of an element. But then $e_v^\lambda$ is actually a graded
automorphism of $R[P]$ because $e_v^{-\lambda}$ is its inverse.

It is clear that $e_v^\lambda$ is just an elementary matrix in the
special case when $P=\Delta_n$, after the identification
$\ggr.aut_R(P)=\GL_{n+1}(R)$. Accordingly, the automorphisms of
type $e_v^\lambda$ are called \emph{elementary}, and the group
they generate in $\ggr.aut_R(P)$ is denoted by $\EE_R(P)$.

\begin{remark}\label{inthisway} Above we have generalized the basic building blocks of higher
$K$-theory of rings to the polytopal setting: general linear
groups and their elementary subgroups. As mentioned in Section
\ref{motivations}, the real motivation for us to pursue the
analogy has been the main result of \cite{BrG1} (Theorem \ref{PLg}
below). It is the polytopal version of the fact that an invertible
matrix over a field can be diagonalized by elementary
transformations on rows (or columns) -- or, putting it in
different words, the group $SK_1$ is trivial for fields.
\end{remark}

\begin{proposition}\label{afemb}
Let $R$ be a ring, $P$ a polytope, and $v_1,\ldots,v_s$ pairwise
different column vectors for $P$ with the same base facet
$F=P_{v_i}$, $i=1,\dots,s$.  Then the mapping
$$
\phi:(R,+)^s\to\ggr.aut_R(P),\qquad
(\lambda_1,\ldots,\lambda_s)\mapsto
e_{v_1}^{\lambda_1}\circ\cdots\circ e_{v_s}^{\lambda_s},
$$
is an embedding of groups. In particular, $e_{v_i}^{\lambda_i}$
and $e_{v_j}^{\lambda_j}$ commute for all $i,j\in\{1,\dots,s\}$,
and the inverse of $e_{v_i}^{\lambda_i}$ is
$e_{v_i}^{-\lambda_i}$.

In the special case, when $R$ is a field the homomorphism $\phi$
is an injective homomorphisms of algebraic groups.
\end{proposition}

For the rest of this subsection we assume that $k$ is a field,
$n=\dim P$, and $\AA(F)$ is the image of the map $\phi$ in
Proposition \ref{afemb}

After $\AA(F)$ we introduce some further subgroups of
$\ggr.aut_k(P)$. First, the $(n+1)$-torus $\TT_{n+1}=(k^*)^{n+1}$
acts naturally on $k[P]$ by restriction of its action on
$k[\ZZ^{n+1}]$ that is given by
$$
(\xi_1,\ldots,\xi_{n+1})(e_i)=\xi_ie_i,\qquad i\in[1,n+1];
$$
here $e_i$ is the $i$-th standard basis vector of $\ZZ^{n+1}$.
This gives rise to an algebraic embedding
$\TT_{n+1}\subset\ggr.aut_k(P)$, and we will identify $\TT_{n+1}$
with its image. It consists precisely of those automorphisms of
$k[P]$ that multiply each monomial by a scalar from $k^*$.

Second, the automorphism group $\Sigma(P)$ of the semigroup $S_P$
is in a natural way a finite subgroup of $\ggr.aut_k(P)$. It is
the group of integral affine transformations mapping $P$ onto
itself.

Third, we have to consider a subgroup of $\Sigma(P)$ defined as
follows. Assume $v$ and $-v$ are both column vectors.  Then for
every point $x\in P\cap\ZZ^n$ there is a unique $y\in P\cap\ZZ^n$
such that $\het_v(x,1)=\het_{-v}(y,1)$ and $x-y$ is parallel to
$v$. The mapping $x\mapsto y$ gives rise to a semigroup
automorphism of $S_P$: it `inverts columns' that are parallel to
$v$. It is easy to see that these automorphisms generate a normal
subgroup of $\Sigma(P)$, which we denote by $\Sigma(P)_{\inv}$.

Finally, $\Col(P)$ is the set of column structures on $P$. Now the
main result of \cite{BrG1} is:

\begin{theorem}\label{PLg}
Let $P$ be an $n$-dimensional polytope and $k$ a field.
\begin{itemize}
\item[(a)]
Every element $\gamma\in\ggr.aut_k(P)$ has a (not uniquely
determined) presentation
$$
\gamma=\alpha_1\circ\alpha_2\circ\cdots\circ\alpha_r\circ\tau\circ\sigma,
$$
where $\sigma\in\Sigma(P)$, $\tau\in\TT_{n+1}$, and
$\alpha_i\in\AA(F_{i})$ such that the facets $F_i$ are pairwise
different and $\#(F_i\cap\ZZ^n)\le \#(F_{i+1}\cap\ZZ^n)$,
$i\in[1,r-1]$.
\item[(b)]
For an infinite field $k$ the connected component of unity
$\ggr.aut_k(P)^0\subset\ggr.aut_k(P)$ is generated by the
subgroups $\AA(F_i)$ and $\TT_{n+1}$. It consists precisely of
those graded automorphisms of $k[P]$ which induce the identity map
on the divisor class group of the normalization of $k[P]$.
\item[(c)] $\dim \ggr.aut_k(P)=\#\Col(P)+n+1$.
\item[(d)]
One has $\ggr.aut_k(P)^0\cap\Sigma(P)=\Sigma(P)_\inv$ and
$$
\ggr.aut_k(P)/\ggr.aut_k(P)^0\approx\Sigma(P)/\Sigma(P)_\inv.
$$
Furthermore, if $k$ is infinite, then $\TT_{n+1}$ is a maximal
torus of $\ggr.aut_k(P)$.
\end{itemize}
\end{theorem}

\section{Stable groups of elementary automorphisms and Polyhedral $K_2$}\label{MILN}

\subsection{Product of column vectors}\label{PRODUCTS}
The product of two column vectors $u,v\in\Col(P)$ is defined as
follows: we say that the product $uv$ exists if $u+v\neq 0$ and
for every point $x\in\L_P\setminus P_u$ the condition $x+u\notin
P_v$ holds. In this case, we define the product as $uv=u+v$. It is
easily seen that $uv\in\Col(P)$ and $P_{uv}=P_u$.

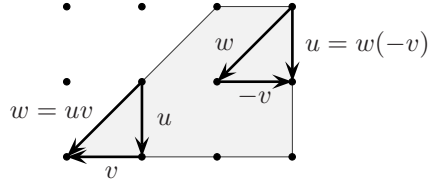
\begin{figure}[htb]
\begin{center}
 \psset{unit=1cm}
 \def\vertex{\pscircle[fillstyle=solid,fillcolor=black]{0.05}}
\begin{pspicture}(-0.3,-0.5)(4,2.5)
 \pspolygon[style=fyp,linecolor=darkgray](0,0)(3,0)(3,2)(2,2)
 \footnotesize
 \multirput(0,1)(1,0){4}{\vertex}
 \multirput(0,0)(1,0){4}{\vertex}
 \multirput(0,2)(1,0){4}{\vertex}
 \psline[style=fatline]{->}(1,1)(0,0)
 \psline[style=fatline]{->}(1,0)(0,0)
 \psline[style=fatline]{->}(2,1)(3,1)
 \psline[style=fatline]{->}(1,1)(1,0)
 \psline[style=fatline]{->}(3,2)(2,1)
 \psline[style=fatline]{->}(3,2)(3,1)
 \rput(-0.2,0.6){$w=uv$}
 \rput(1.3,0.5){$u$}
 \rput(0.6,-0.2){$v$}
 \rput(2.5,0.8){$-v$}
 \rput(2.1,1.5){$w$}
 \rput(4.0,1.5){$u=w(-v)$}
\end{pspicture}
\end{center}
\caption{The product of two column vectors} \label{ProdCol}
\end{figure}

Figure \ref{ProdCol} shows a polytope with all its column vectors
and the two existing products $w=uv$ and $u=w(-v)$.

In the case of a unimodular simplex the product of two oriented
edges, viewed as column vectors, exists if and only if they are
not opposite to each other and the end point of the first edge is
the initial point of the second edge.

\subsection{Balanced polytopes}\label{BALANCED}

A polytope $P$ is called \emph{balanced} if $\la P_u,v\ra\leq 1$
for all $u,v\in\Col(P)$. One easily observes that $P$ is balanced
if and only if $|\la P_u,v\ra|\leq 1$ for all $u,v\in\Col(P)$.

The reason we introduce balanced polytopes is that the main
results of \cite{BrG5,BrG6} are only proved for this class of
polytopes. However, it is not yet excluded that everything
generalizes to arbitrary polytopes.

We give the classification result in dimension 2. It uses the
notion of projective equivalence: $n$-dimensional polytopes
$P,Q\subset\RR^n$ are called \emph{projectively equivalent} if and
only if $P$ and $Q$ have the same dimension, the same
combinatorial type, and the faces of $P$ are parallel translates
of the corresponding ones of $Q$. An alternative definition in
terms of \emph{normal fans} is given in Subsection \ref{FUNCT}.

Recall the notation
$\Delta_n=\conv(O,(1,\ldots,0),\ldots,(0,\ldots,1))$ for the unit
$n$-simplex.

\begin{theorem}\label{DIM2}
For a balanced polygon $P$ there are exactly the following
possibilities (up to integral-affine equivalence):
\begin{itemize}
\item[(a)] $P$ is a multiple of the unimodular triangle
$P_a=\Delta_2$. Hence $\Col(P)= \{\pm u,\pm v,\allowbreak\pm w\}$
and the column vectors are subject to the obvious relations,
\item[(b)] $P$ is projectively equivalent to the trapezoid
$P_b=\conv\bigl((0,0),(2,0),(1,1),\allowbreak(0,1)\bigr)$, hence
$\Col(P)=\{u,\pm v,w\}$ and the relations in $\Col(P)$ are $uv=w$
and $w(-v)=u$, \item[(c)] $\Col(P)=\{u,v,w\}$ and $uv=w$ is the
only relation, \item[(d)] $\Col(P)$ has any prescribed number of
column vectors, they all have the same base edge (clearly, there
are no relations between them), \item[(e)] $P$ is projectively
equivalent to the unit lattice square $P_e$, hence $\Col(P)=\{\pm
u,\pm v\}$ with no relations between the column vectors,
\item[(f)] $\Col(P)=\{u,v\}$ so that $P_u\neq P_v$ with no
relations in $\Col(P)$.
\end{itemize}
\end{theorem}

It turns out that polyhedral $K$-groups are invariants of the
projective equivalence classes of polytopes (in arbitrary
dimension); see Proposition \ref{PrEq} below.

\subsection{Doubling along a facet}\label{DOUBLING}

Let $P\subset\RR^n$ be a polytope and $F\subset P$ be a facet. For
simplicity we assume that $0\in F$, a condition that can be
satisfied by a parallel translation of $P$. Denote by
$H\subset\RR^{n+1}$ the $n$-dimensional linear subspace that
contains $F$ and whose normal vector is perpendicular to that of
$\RR^n=\RR^n\oplus0\subset\RR^{n+1}$ (with respect to the standard
scalar product on $\RR^{n+1}$). Then the upper half space
$H\cap\bigl(\RR^n\times\RR_+\bigr)$ contains a congruent copy of
$P$ which differs from $P$ by a $90^\circ$ rotation. Denote the
copy by $P^{\vt_F}$, or just by $P^\vt $ if there is no danger of
confusion.

Note that $P^\vt$ is not always a lattice polytope with respect to
the standard lattice $\ZZ^{n+1}$. However, it is so with respect
to the sublattice $(\ZZ^n)^{\vt_F}$ which is the image of $\ZZ^n$
under the $90^\circ$ rotation.

The operator of doubling along a facet is then defined by
$$
P^{\sq_F}=\conv(P,P^\vt )\subset \RR^{n+1}.
$$

The doubled polytope is a lattice polytope with respect to the
subgroup $(\ZZ^n)^{\sq_F}=\ZZ^n+(\ZZ^n)^{\vt_F}\subset\RR^{n+1}$.
After a change of basis in $\RR^{n+1}$ that does not affect
$\RR^n$ we can replace $(\ZZ^n)^{\sq_F}$ by $\ZZ^{n+1}$, and
consider $P^{\sq_F}$ as an ordinary lattice polytope in
$\RR^{n+1}$. In what follows, whenever we double a lattice
polytope $P\subset\RR^n$ along a facet $F$, the lattice of
reference in $\RR^{n+1}$ is always $\ZZ^n+(\ZZ^n)^{\vt_F}$. For
simplicity of notation this lattice will be denoted by
$\ZZ^{n+1}$.


\begin{figure}[htb]
\begin{center}
\psset{unit=1cm}
\def\vertex{\pscircle[fillstyle=solid,fillcolor=black]{0.07}}
\begin{pspicture}(-2.8,0)(1,2)
\psset{viewpoint=3.5 2 -1.5} \footnotesize \ThreeDput[normal=0 0
-1](0,0,0){
  \pspolygon[style=fyp](0,0)(3,0)(2,1.5)(0,1.5)(0,0)
  \multirput(0,0)(1,0){4}{\vertex}
  \multirput(0,1.5)(1,0){3}{\vertex}
  \psline[linewidth=1.5pt]{->}(0,1.5)(0,0) 
 }
\ThreeDput[normal=0 -1 0](0,0,0){
  \pspolygon[style=fyp](0,0)(3,0)(2,1.5)(0,1.5)(0,0)
  \multirput(0,1.5)(1,0){3}{\vertex}
  \psline[linewidth=2.0pt]{->}(0,1.5)(0,0) 
 }
\ThreeDput[normal=-1 0 0](0,0,0){
  \psline[linewidth=1.0pt]{->}(1.5,0)(0,1.5)
 }
\ThreeDput[normal=-1 0 0](2,0,0){
  \psline[linewidth=1.0pt]{->}(0,1.5)(1.5,0)
 }
 \rput(-2.2,-0.0){$P$}
 \rput(-0.25,2.0){$P^\vt $}
 \rput(-1.3,0.5){$F$}
 \rput(-0.6,-0.25){$v^-$}
 \rput(0.3,0.5){$v^\vt $}
\end{pspicture}
\caption{Doubling along the facet $F$} \label{Doubling}
\end{center}
\end{figure}
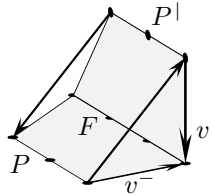

In case $F=P_v$ for some $v\in\Col(P)$ we will use the notation
$P^{\sq_F}=P^{\sq_v}$.

\subsection{The stable group of elementary automorphisms}\label{ELAUT}
An ascending infinite chain of lattice polytopes ${\mathfrak
P}=(P=P_0\subset P_1\subset\dots)$ is called a \emph{doubling
spectrum} if the following conditions hold:
\begin{itemize}
\item[(i)]
for every $i\in\ZZ_+$ there exists a column vector
$v\subset\Col(P_i)$ such that $P_{i+1}=P_i^{\sq_v}$,
\item[(ii)]
for every $i\in\ZZ_+$ and any $v\in\Col(P_i)$ there is an index
$j\geq i$ such that $P_{j+1}=P_j^{\sq_v}$.
\end{itemize}

Here we use the natural inclusion $\Col(P_i)\subset\Col(P_{i+1})$.

One says that $v\in\Col(P_i)$ is \emph{decomposed} at the $j$th
step in ${\frak  P}$ for some $j\geq i$ if $P_{j+1}=P_j^{\sq_v}$.

Associated to a doubling spectrum $\frak P$ is the `infinite
polytopal' algebra
$$
R[\fP]=\lim_{i\to\infty}R[P_i]
$$
and the filtered union
$$
\Col(\fP)=\lim_{i\to\infty}\Col(P_i).
$$
The product of two vectors from $\Col(\fP)$ is defined in the
obvious way, using the definition for a single polytope. Also, we
can speak of systems of elements of $\Col(\fP)$ having the same
base facets, etc.

Elements $v\in\Col(\fP)$ and $\lambda\in R$ give rise to a graded
automorphism of $R[\fP]$ as follows: we choose an index $i$ big
enough so that $v\in\Col(P_i)$. Then the elementary automorphisms
$e_v^{\lambda}\in\EE_R(P_j)$, $j\geq i$ constitute a compatible
system and, therefore, define a graded automorphism of $R[{\frak
P}]$. This automorphism will also be called `elementary' and it
will be denoted by $e_v^{\lambda}$.

The group $\EE(R,\fP)$ is by definition the subgroup of
$\ggr.aut_R(R[\fP])$, generated by all elementary automorphisms.

\begin{remark}\label{nolimit}
Unlike the classical situation of unimodular simplices,
the group $\EE(R,\fP)$ can \emph{not} be represented
as a direct limit of the `unstable' groups $\EE_R(P_i)$,
$i\in\ZZ_+$.
\end{remark}

\begin{theorem}\label{elaut}
Let $R$ be a ring and $P$ be a polytope (not necessarily balanced)
admitting a column structure. Assume $\fP=(P\subset P_1\subset
P_2\subset\dots)$ is a doubling spectrum. Then:
\begin{itemize}
\item[(a)] $\EE(R,\fP)$ is naturally isomorphic to
$\EE(R,{\mathfrak Q})$ for any other doubling spectrum ${\mathfrak
Q}=(P\subset Q_1\subset Q_2\subset\dots)$.

\item[(b)] $\EE(R,\fP)$ is perfect.

\item[(c)] The center of $\EE(R,\fP)$ is trivial.

\item[(d)] $e_u^{\lambda}\circ e_u^{\mu}=e_u^{\lambda+\mu}$ for
every $u\in\Col(\fP)$ and $\lambda,\mu\in R$. \item[(e)] If $P$ is
balanced, $u,v\in\Col(\fP)$, $u+v\neq 0$ and $\lambda,\mu\in R$
then
$$
[e_u^\lambda,e_v^\mu]=
\begin{cases}e_{uv}^{-\lambda\mu}&\text{if}\ uv\ \text{exists},\\
1&\text{if}\ u+v\notin\Col(\fP).
\end{cases}
$$
\end{itemize}
\end{theorem}

The difficult parts of this theorem are the claims (c) and (e),
which in the special case $P=\Delta_n$ are just standard facts.

Thanks to Theorem \ref{elaut}(a) we can use the notation
$\EE(R,P)$ for $\EE(R,\fP)$.

\begin{remark}\label{steinberg}
Theorem \ref{elaut}(e) is the generalization of Steinberg's
relations between elementary matrices to balanced polytopes.
\end{remark}

\subsection{The Schur multiplier}\label{SCHUR}

Let $P$ be a balanced polytope and $\fP=(P\subset P_1 \subset
P_2\subset\dots)$ be a doubling spectrum. Then for a ring $R$ we
define the \emph{stable polytopal Steinberg group} $\St(R,P)$ as
the group generated by symbols $x_v^\lambda$, $v\in\Col({\mathfrak
P})$, $\lambda\in R$, which are subject to the relations
$$
x_v^\lambda x_v^\mu=x_v^{\lambda+\mu}
$$
and
$$
[x_u^\lambda,x_v^\mu]=
\begin{cases}
x_{uv}^{-\lambda\mu}&\text{if}\ uv\ \text{exists,}\\
\\
1&\text{if}\ u+v\notin\Col(\fP)\cup\{0\}.
\end{cases}
$$
The use of the notation $\St(R,P)$ is justified by the fact that,
like in Theorem \ref{elaut}(a), the stable Steinberg groups are
determined by the underlying doubling spectra (with the same
initial polytope) up to canonical isomorphism.

The central result of \cite{BrG5} is the following

\begin{theorem}\label{schur}
For a ring $R$ and a balanced polytope $P$ the natural surjective
group homomorphism $\St(R,P)\to\EE(R,P)$ is a universal central
extension whose kernel coincides with the center of $\St(R,P)$.
\end{theorem}

The group $\Ker\bigl(\St(R,P)\to\EE(R,P)\bigr)$ is called the
\emph{polyhedral Milnor group}. We denote it by $K_2(R,P)$.
Clearly, when $P$ is a unimodular simplex $K_2(R,P)$ is the usual
Milnor group $K_2(R)$ \cite{Mi}.

\section{Rigid systems of column vectors}\label{RIGID}

We can speak of the product $\prod_{i=1}^mv_i$ of elements
$v_i\in\Col(P)$ whenever the following two conditions are
satisfied:
\begin{itemize}
\item[(i)] the products $v_iv_{i+1}$ exist for all $i\in[1,m-1]$,

\item[(ii)] $\sum_{i=r}^sv_i\neq 0$ for all $1\leq r<s\leq m$.
\end{itemize}
In this case every bracketing of the sequence $v_1v_2\dots v_m$
yields pairs of column vectors whose products exist.

It is useful to have another, weaker notion of product. We say
that $\prod_{i=1}^mv_i$ exists \emph{weakly} if there is a
bracketing of the sequence
$$
v_1v_2\cdots v_m
$$
such that all the recursively defined products of pairs of column
vectors exist. Since $v_1\cdots v_n=v_1+\dots + v_n$ in the case
of weak existence, the value of the product does not depend on the
bracketing.

By $\la V\ra$ we denote the hull of $V$ in $\Col(P)$ under
products (of two column vectors). One has $v\in\la V\ra$ if and
only if there exist $v_1,\dots,v_m\in V$ such that $v=v_1\cdots
v_m$ is their weak product.

For simplicity we introduce the following convention: $v_1\cdots
v_m\in [V]$ means that the product of $v_1,\dots, v_m$ exists (in
the strong sense), whereas $v_1\cdots v_m\in \la V\ra$ means that
the product of $v_1,\dots, v_m$ exists in the weak sense.

We will represent certain partial product structures on sets of
column vectors by equivalence classes of directed paths in graphs.
The \emph{graphs} considered here are finite directed graphs $\bG$
satisfying the following conditions:
\begin{itemize}
\item[(i)] $\bG$ has no isolated vertices;
\item[(ii)] $\bG$ has no multiple edges and no edges from a vertex
to itself;
\item[(iii)] if vertices $a$ and $b$ are connected by an edge, then
there is no other directed path connecting $a$ and $b$.
\end{itemize}
Condition (iii) implies that there are no directed cycles in $\bG$
(but the existence of non-directed cycles is not excluded). A
\emph{path} is always assumed to be oriented.

The set of nonempty paths in a graph $\bF$ carries a natural
partial product structure -- $ll'$ exists if the end point of the
path $l$ is the initial point for $l'$. The set of all paths in
$\bF$ is denoted by $\path{\bF}$. There is an equivalence relation
on $\path{\bf F}$: two paths are considered to be equivalent if
they have the same initial and the same end point. We let
$\overline{\path{\bF}}$ denote the corresponding quotient set.

\begin{definition}\label{rigid}
A system of column vectors $V\subset\Col(P)$ is called
\emph{rigid} if the following conditions are satisfied:
\begin{itemize}
\item[(a)] $[V]$ does not contain a subset of type $\{v,-v\}$,
$v\in\Col(P)$;

\item[(b)] $[V]=\la V\ra$;

\item[(c)] there exist a graph $\bF$ and an isomorphism
$[V]\approx\overline{\path{\bF}}$ of partial product structures.
\end{itemize}
\end{definition}

\section{Higher polyhedral $K$-groups}\label{HPK}

In this section we assume that $R$ is a ring and $P$ is a balanced
polytope admitting a column structure.

\subsection{Triangular subgroups in $\EE(R,P)$ and
$\St(R,P)$}\label{TRIANG}

We fix a doubling spectrum $\fP=(P\subset P_1 \subset\cdots)$.
Thanks to Theorem \ref{elaut}(a) (and its straightforward analogue
for polyhedral Steinberg groups) all the objects defined below are
independent of the fixed spectrum.

We say that $V\subset\Col(\fP)$ is a rigid system if there exists
an index $j\in\NN$ such that $V$ is a subset of $\Col(P_j)$ and is
rigid.

\begin{definition}\label{triang}\leavevmode\par
\begin{itemize}
\item[(a)] A subgroup $G\subset\EE(R,P)$ is called
\emph{triangular} if there exists a rigid system
$V\subset\Col(\fP)$ such that $G$ is generated by the elementary
automorphisms $e_v^\lambda$, where $\lambda$ runs through $R$ and
$v$ through $V$. The triangular subgroup corresponding to a rigid
system $V$ is denoted by $G(R,V)$, and $\T(R,P)$ is the family of
all triangular subgroups of $\EE(R,P)$.

\item[(b)] The triangular subgroups of $\St(R,P)$ are defined
similarly.
\end{itemize}
\end{definition}

\subsection{Volodin's theory}\label{VOLOD}

\begin{definition}\label{volod}
\leavevmode\par
\begin{itemize}
\item[(a)]
The $d$-simplices of the \emph{Volodin simplicial set}
$\vv(\EE(R,P))$ are those sequences
$(\epsilon_0,\dots,\epsilon_d)\in(\EE(R,P))^{d+1}$ for which there
exists a triangular group $G\in\T(R,P)$  such that
$\epsilon_k\epsilon_l^{-1}\in G$, $k,l\in[0,d]$. The $i$th face
(resp.\ degeneracy) of $\vv(\EE(R,P))$ is obtained by omitting
(resp.\ repeating) $\epsilon_i$.
\item[(b)]
The simplicial set $\vv(\St(R,P))$ is defined analogously.
\item[(c)]
The higher Volodin polyhedral $K$-groups of $R$ are defined by
$$
K_i^{\V}(R,P)=\pi_{i-1}\bigl(|\vv(\EE(R,P))|,({\text{\bf
Id}})\bigr),\quad i\geq2.
$$
\end{itemize}
where $|-|$ refers to the geometric realization of a simplicial
set.
\end{definition}

The definition of the Volodin simplicial set is independent of the
choice of $\frak P$ and one has
$$
K_i^{\V}(R,P)=\pi_{i-1}(\vv(\St(R,P))),\quad i\geq3.
$$
When $P$ is a unimodular simplex of arbitrary dimension Definition
\ref{volod} gives the usual Volodin theory \cite{Vo}.

\subsection{Quillen's theory}\label{QUIL}

We define \emph{Quillen's higher polyhedral $K$-groups} by
$$
K_i^{\Qu}(R,P)=\pi_i(\B\EE(R,P)^+),\qquad i\geq2,
$$
where $\B\EE(R,P)^+$ refers to Quillen's $+$ construction applied
to $\B\EE(R,P)$ with respect to the whole group
$\EE(R,P)=[\EE(R,P),\EE(R,P)]$ (Theorem \ref{elaut}(b)).

We have the equalities
$$
K_i^{\Qu}(R,P)=\pi_i(\B\St(R,P)^+),\qquad i\geq3,\label{QuilEq}
$$
where the $+$ construction is considered with respect to the whole
group $\St(R,P)$.

\begin{proposition}\label{v=q2}
$K_2^{\Qu}(R,P)=K_2(R,P)=K^{\V}_2(R,P)$.
\end{proposition}

For a unimodular simplex $P=\Delta_n$ we recover Quillen's theory
\cite{Qu1}.

\subsection{Functorial properties}\label{FUNCT}

Let $Q$ be another balanced polytope. If there exists a mapping
$\mu:\Col(P)\to\Col(Q)$, such that the conditions
$$
\text{(i)}\quad\la P_w,v\ra=\la
Q_{\mu(w)},\mu(v)\ra\qquad\text{and}\qquad\text{(ii)}\quad
\mu(vw)=\mu(v)\mu(w)\text{ if $vw$ exists,}
$$
hold for all $v,w\in\Col(P)$, then the assignment
$x_v^\lambda\mapsto x_{\mu(v)}^\lambda$ induces a homomorphism
$$
\St(R,\mu):\St(R,P)\to\St(R,Q).
$$
Moreover, if $\mu$ is bijective, then
$$
\St(R,P)\approx\St(R,Q),\quad\EE(R,P)\approx\EE(R,Q),\quad
K_2(R,P)\approx K_2(R,Q).
$$
This observation allows one to study polyhedral $K$-theory as a
functor also in the polytopal argument. The map $\mu$ is called a
\emph{$K$-theoretic morphism} from $P$ to $Q$. Though we cannot
prove $K_2$-functoriality for \emph{all} maps $\mu$, it is useful
to note the $\St$-functoriality, since it implies bifunctoriality
of the higher polyhedral $K$-groups with covariant arguments:
\begin{multline*}
K_i^{\Qu}(-,-),K_i^{\V}(-,-):\underline{\text{\it Commutative
Rings}}
\ \times\ \underline{\text{\it Balanced Polytopes}}\to\\
\to\underline{\text{\it Abelian Groups}},\quad i\geq3.
\end{multline*}

The \emph{normal fan} ${\cal N}(P)$  of a finite convex (not
necessarily lattice) polytope $P\subset\RR^n$ is defined as the
complete fan in the dual space $(\RR^n)^*=\Hom(\RR^n,\RR)$ given
by the system of cones
$$
\bigl(\{\phi\in(\RR^n)^*\mid\max_P(\phi)=F\},\ F\ \text{a face
of}\ P\bigr).
$$
Two polytopes $P,Q\subset\RR^n$ are projectively equivalent (see
Section \ref{BALANCED}) if and only if ${\cal N}(P)={\cal N}(Q)$.

\begin{proposition}\label{PrEq}
If $P$ and $Q$ are projectively equivalent balanced polytopes,
then $K_i^{\Qu}(R,P)\approx K_i^{\Qu}(R,Q)$ and
$K_i^{\V}(R,P)\approx K_i^{\V}(R,Q)$ for $i\geq2$.
\end{proposition}

\section{On the coincidence of Quillen's and Volodin's theories}\label{Q=V}

All polytopes are assumed to be balanced and to admit a column
vector, unless specified otherwise.

\begin{figure}[htb]
\begin{center}
\footnotesize
 \psset{unit=1cm}
\begin{pspicture}(-1,0)(1,2)
 \psline{->}(-1,0)(0,1)
 \psline{->}(1,0)(0,1)
 \psline{->}(0,1)(0,2)
 \psline[linestyle=dashed]{->}(-1,0)(1,0)
 \rput(-0.6,0.7){$a$}
 \rput(0.6,0.7){$b$}
 \rput(-0.2,1.5){$c$}
 \rput(0,0.25){$d$}
\end{pspicture}
\quad
\begin{pspicture}(-1,0)(1,2)
 \psline{->}(-1,0)(0,1)
 \psline{->}(1,0)(0,1)
 \psline{->}(0,1)(0,2)
 \psline[linestyle=dashed]{->}(1,0)(-1,0)
 \rput(-0.6,0.7){$a$}
 \rput(0.6,0.7){$b$}
 \rput(-0.2,1.5){$c$}
 \rput(0,0.25){$d$}
\end{pspicture}
\quad
\begin{pspicture}(-1,0)(1,2)
 \psline{->}(0,0)(-1,1)
 \psline{->}(0,0)(1,1)
 \psline{->}(-1,1)(0,2)
 \psline{->}(1,1)(0,2)
 \psline[linestyle=dashed](-1,1)(1,1)
 \psline[linestyle=dashed]{->}(-0.05,1)(0.05,1)
 \rput(-0.6,0.4){$a$}
 \rput(0.6,0.4){$c$}
 \rput(-0.6,1.65){$b$}
 \rput(0.6,1.65){$d$}
 \rput(0,1.25){$t$}
\end{pspicture}
\quad
\begin{pspicture}(-1,0)(1,2)
 \psline{->}(0,0)(-1,1)
 \psline{->}(0,0)(1,1)
 \psline{->}(-1,1)(0,2)
 \psline{->}(1,1)(0,2)
 \psline[linestyle=dashed](1,1)(-1,1)
 \psline[linestyle=dashed]{->}(0.05,1)(-0.05,1)
 \rput(-0.6,0.4){$a$}
\rput(0.6,0.4){$c$}
 \rput(-0.6,1.65){$b$}
 \rput(0.6,1.65){$d$}
 \rput(0,1.25){$t$}
\end{pspicture}
\end{center}
\caption{$\Col$-divisibility} \label{FigColDiv}
\end{figure}
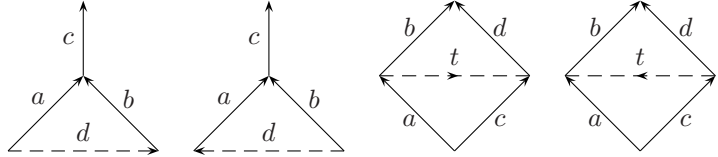

\begin{definition}\label{divis}
A (balanced) polytope $P$ is \emph{$\Col$-divisible} if its column
vectors satisfy the following condition:
\begin{itemize}
\item[(\cd1)]
if $ac$ and $bc$ exist and $a\neq b$, then $a=db$ or $b=da$ for
some $d$;
\item[(\cd2)]
if $ab=cd$ and $a\not=c$, then there exists $t$ such that $at=c$,
$td=b$, or $ct=a$, $tb=d$.
\end{itemize}
(See Figure \ref{FigColDiv}.)
\end{definition}

The main result of \cite{BrG6} is the following

\begin{theorem}\label{q=v}
Suppose $P$ is a $\Col$-divisible polytope. Then
$$
K_i^{\Qu}(R,P)=K_i^{\V}(R,P),\qquad i\geq2.
$$
\end{theorem}

The proof is a `polytopal extension' of Suslin's proof \cite{Su}
of the coincidence of the usual theories.

However, we expect that Quillen's and Volodin's theories diverge
for general balanced polytopes, see Remark \ref{final}.

\section{Computations}\label{POLYG}

\subsection{The case of polygons} The class of $\Col$-divisible
polytopes may at first glance seem rather restricted. However, it
follows immediately from Theorem \ref{DIM2} that \emph{all}
balanced polytopes of dimension 2 are $\Col$-divisible.

Let $R$ be a ring. In Theorem \ref{DIM2} we have grouped all
balanced polygons in six infinite series which give rise to the
following isomorphism classes of stable elementary automorphism
groups:
\begin{equation}
\tag{a}\EE_a=\E(R),
\end{equation}
\begin{equation}
\tag{b}\EE_b=
\begin{pmatrix}
\E(R)&\End_R(\oplus_{\NN}R)\\
&\\
0&\E(R)
\end{pmatrix},
\end{equation}
\begin{equation}
\tag{c}\EE_c=
\begin{pmatrix}
\E(R)&\End_R(\oplus_{\NN}R)&
\Hom_R(\oplus_{\NN}R,R)\\
&\\
0&\E(R)&\Hom_R(\oplus_{\NN}R,R)\\
&\\
0&0&1
\end{pmatrix},
\end{equation}
\begin{equation}
\tag{d} \EE_{d,t}=
\begin{pmatrix}
\E(R)&\Hom_R(\oplus_{\NN}R,R^t)\\
0&\text{\bf Id}_t
\end{pmatrix},\quad t\in\NN,
\end{equation}
\begin{equation}
\tag{e}\EE_e=\E(R)\times\E(R),
\end{equation}
\begin{equation}
\tag{f}\EE_f=
\begin{pmatrix}
\E(R)&\Hom_R(\oplus_{\NN}R,R)\\
&\\
0&1
\end{pmatrix}
\times
\begin{pmatrix}
\E(R)&\Hom_R(\oplus_{\NN}R,R)\\
&\\
0&1
\end{pmatrix}.
\end{equation}

\begin{definition}
A ring $R$ is an \emph{$S(n)$-ring} if there are $r_1,\dots,r_n\in
R^*$ such that the sum of each nonempty subfamily is a unit. If
$R$ is an $S(n)$-ring for all $n\in\NN$, then $R$ has \emph{many
units}.
\end{definition}

The class of rings with many units includes local rings with
infinite residue fields and algebras over rings with many units.

\begin{theorem}\label{comput}
For every ring $R$ and every index $i\geq2$ we have:
\begin{itemize}
\item[(a)]
$\pi_i(\B\EE_a^+)=K_i(R)$,
\item[(b)]
$\pi_i(\B\EE_b^+)=K_i(R)\oplus K_i(R)$,
\item[(c)]
$\pi_i(\B\EE_c^+)=K_i(R)\oplus K_i(R)$ if $R$ has many units,
\item[(d)]
$\pi_i(\B\EE_{d,t}^+)=K_i(R)$ if $R$ has many units,
\item[(e)]
$\pi_i(\B\EE_e^+)=K_i(R)\oplus K_i(R)$,
\item[(f)]
$\pi_i(\B\EE_f^+)=K_i(R)\oplus K_i(R)$ if $R$ has many units.
\end{itemize}
\end{theorem}

The proof is based on homological computations for the
corresponding matrix groups due to
Nesterenko-Suslin \cite{NSu} and Quillen \cite{Qu2}.

\subsection{Higher dimensional polytopes} It seems that a similar `almost triangular' matrix group
interpretation is possible for the group of elementary
automorphisms for \emph{all} $\Col$-divisible polytopes. Then,
based on the techniques of Berrick and Keating \cite{BKe,Ke}, the
corresponding $K$-groups should be computable in terms of the
usual $K$-groups of the underlying ring. This remark leads us to
the following

\begin{conjecture}\label{conjecture}
For a commutative ring $R$ and a $\Col$-divisible polytope $P$ of
arbitrary dimension we have
$$
K_i(R,P)=\underbrace{K_i(R)\oplus\cdots\oplus K_i(R)}_{{\mathfrak
c}(P)},\qquad i\ge2,
$$
where ${\mathfrak c}(P)\le\dim P$ is a natural number explicitly
computable in terms of the partial product table of $\Col(P)$.
\end{conjecture}

\begin{remark}\label{final}
For balanced but not $\Col$-divisible polytopes we may expect that
Quillen's and Volodin's theories diverge and we get really new
$K$-groups. The simplest candidate for such a deviation from the
usual theory is the pyramid over the unit square shown below --
its column vectors are the four oriented edges of the square and
four oriented edges emerging from the top vertex. This polytope
has shown up several times in our papers as a counterexample to
several natural conditions.
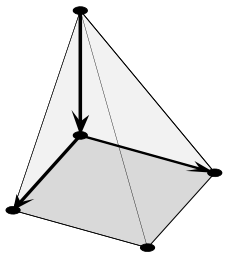
\begin{figure}[htb]
\begin{center}
\psset{unit=2cm}
\def\vertex{\pscircle[fillstyle=solid,fillcolor=black]{0.05}}
\begin{pspicture}(-0.5,-0.7)(1,0.7)
\psset{viewpoint=-1 -2 -1.5} \ThreeDput[normal=0 0 1](0,0,0){
  \pspolygon[linewidth=0pt, style=fyp, fillcolor=medium](1,0)(0,0)(0,1)(1,1)
} \ThreeDput[normal=0 -1 0]{
  \pspolygon[linewidth=0pt, style=fyp](1,0)(0,0)(0,1)
} \ThreeDput[normal=1 0 0]{
  \pspolygon[linewidth=0pt, style=fyp](1,0)(0,0)(0,1)
} \ThreeDput[normal=0 0 1](0,0,0){
  \rput(0,0){\vertex}
  \rput(0,1){\vertex}
  \rput(1,0){\vertex}
  \rput(1,1){\vertex}
  \psline[linestyle=dashed](0,1)(0,0)
  \psline[linewidth=1.5pt]{->}(0,0)(1,0)
  \psline[linewidth=1.5pt]{->}(0,0)(0,1)
  \psline(0,1)(1,1)
} \ThreeDput[normal=0 -1 0]{
  \psline(0,1)(1,0)
  \psline[linewidth=1.5pt]{->}(0,1)(0,0)
}
 \ThreeDput[normal=1 0 0]{\psline(0,1)(1,0)}
 \ThreeDput[normal=1 -1 0]{\psline(0,1)(1.41,0)}
 \ThreeDput(0,0,1){\rput(0,0){\vertex}}
\end{pspicture}
\end{center}
\caption{The pyramid over the unit square} \label{ThreeCol}
\end{figure}
\end{remark}

\end{document}